\documentclass{aero}

\usepackage[utf8]{inputenc}
\usepackage{textcomp}

\usepackage{graphicx}
\usepackage{amsmath}
\usepackage[version=4]{mhchem}
\usepackage{siunitx}
\usepackage{longtable,tabularx}
\usepackage[boxed, ruled]{algorithm2e}
\usepackage{mathtools}
\usepackage{arydshln}

\usepackage{multirow}
\usepackage{setspace}

\DeclareMathOperator*{\Minimize}{Minimize}
\DeclareMathOperator*{\diag}{diag}

\ADsetup{
title    = {Low-Thrust Trajectory Optimization with Quantum Computing and Sequential Convex Programming},
author   = {Carmine Giordano$^{1}$\cor{}},
address  = {
  1.  Dept.\ of Aerospace Science and Technology, Politecnico di Milano, Milano, Italy, 20156
},
email    = {carmine.giordano@polimi.it},
abstract = {
  Low-thrust trajectory optimization is a central task in interplanetary mission design, but its nonlinear dynamics and operational constraints often lead to challenging non-convex optimal-control problems. Sequential convex programming has emerged as an effective approach to address these problems, while quantum annealing offers a complementary paradigm for solving quadratic unconstrained binary optimization problems. This paper introduces quSCP, a quantum-based sequential convex programming framework that reformulates each convex subproblem as a quadratic unconstrained binary optimization problem suitable for quantum and hybrid quantum--classical solvers. Equality, inequality, and trust-region constraints are embedded through quadratic penalty terms, while an iterative refinement strategy is used to reduce the accuracy loss introduced by binary discretization. The method is assessed on a fuel-optimal Earth--Mars low-thrust transfer by comparing standard sequential convex programming, a continuous quadratic unconstrained formulation, direct quantum processing unit sampling, and D-Wave hybrid solvers. Results show that quSCP produces physically consistent trajectories with propellant consumption and nonlinear constraint violations close to classical benchmarks. Direct quantum annealing is feasible only for small instances because of embedding overhead and hardware connectivity limits, whereas hybrid solvers scale to larger discretizations. Although no computational quantum advantage is demonstrated with current hardware, the results show that quantum and hybrid quantum--classical optimization can already provide competitive solutions for demanding trajectory design problems.
  },
keywords = {
  Trajectory Optimization \sep Quantum Computing \sep Quantum Annealing \sep Sequential Convex Programming
},
}

\begin{document}

\maketitle  

\section{Introduction}
The consumption, time of flight, or maximize mission robustness is of capital importance for both robotic and human missions beyond Earth. In recent years, the rapid advancement of computational capabilities has enabled the development of new tools and methodologies aimed at facilitating the implementation and solution of complex trajectory optimization problems in space. The combination of novel computational techniques with increasingly powerful hardware platforms has significantly enhanced the efficiency and flexibility of trajectory design. Methods leveraging artificial intelligence \cite{izzo2019survey,habib2022artificial}, machine learning \cite{zavoli2021reinforcement,federici2021deep}, and general-purpose graphics processing units \cite{arora2015parallel,antony2017rapid,masat2023gpu,chari2024fast} have opened new horizons for solving large-scale and computationally intensive optimization tasks. Among these emerging methodologies, Sequential Convex Programming (SCP) has gained prominence as an effective strategy for addressing non-convex trajectory optimization problems \cite{sagliano2018pseudospectral,hofmann2021rapid,malyuta2022convex,giordano2025goal}. SCP iteratively linearizes the problem around a reference trajectory, solving a sequence of convex subproblems that converge toward an optimal solution. Its application in deep-space trajectory optimization has demonstrated promising results, particularly for low-thrust transfers where the trajectory dynamics are inherently nonlinear and complex \cite{hofmann2021rapid}. Parallel to these algorithmic advances, quantum computing has emerged as a disruptive technology with the potential to tackle computational problems considered intractable for classical computers \cite{nielsen2010quantum}. In particular, Quantum Annealing (QA) \cite{farhi2000quantum} offers an alternative paradigm for solving hard combinatorial optimization problems by mapping them into a Quadratic Unconstrained Binary Optimization (QUBO) form \cite{kadowaki1998quantum}. Although still in its infancy, quantum computing has shown promise in a variety of contexts \cite{peng2008quantum,montanaro2015quantum,zanger2021quantum,yarkoni2022quantum}, in general optimization applications \cite{shukla2019trajectory,hadfield2019quantum,abbas2024challenges}, and it is now being explored within the aerospace community for trajectory planning \cite{makarov2024quantum} and optimization applications \cite{degrossi2025transcription}.\\
In this work, a novel approach that combines Sequential Convex Programming with QUBO-based optimization, labeled as quSCP, is presented. The objective of this study is not to demonstrate a computational quantum advantage with current quantum-annealing hardware, but rather to assess whether low-thrust trajectory optimization subproblems can be reformulated and solved within a quantum-compatible optimization framework. In this sense, the proposed method should be interpreted as a feasibility study and as a methodological step toward future quantum-assisted trajectory design tools. The proposed methodology is applied to a low-thrust Earth-to-Mars transfer scenario, serving as a benchmark to evaluate the quality of the solutions obtained with QPU-based and hybrid quantum--classical solvers. Although current devices impose severe limitations in terms of embeddable problem size, coefficient precision, and overall runtime, the ability to obtain physically meaningful solutions from a QUBO formulation is relevant for future applications in which quantum optimization may become more effective. This includes higher-level mission-design problems, such as low-thrust multiple-gravity-assist trajectory design or campaign-level optimization, where the low-thrust transfer optimization considered here may appear as one component of a larger discrete-continuous quantum optimization pipeline.\\
The paper is organized as follows. First, the background about SCP and QA is presented in Section \ref{sec:Background}. Then, Section \ref{sec:QUBO} shows how the low-thrust optimization problem is adapted to the QUBO framework. The test-case scenario, results, and comparisons are shown in Section \ref{sec:Results}. Section \ref{sec:Conclusions} concludes the work.

\section{Background}\label{sec:Background}

\subsection{Sequential Convex Programming}\label{subsec:SCP}
Sequential convex programming is a powerful method for solving optimal control problem by transforming nonconvex problems into a series of easier-to-solve convex subproblems \cite{malyuta2022convex}. This approach allows efficient handling of complex dynamics and constraints, making it particularly useful for optimizing continuous-thrust trajectories. The efficiency and flexibility of SCP have made this approach popular in various aerospace contexts, such as powered descent guidance \cite{acikmese2007convex} and interplanetary trajectory design \cite{hofmann2021rapid, hofmann2023performance}.\\ 
The generic optimal control problem for a spacecraft equipped with low-thrust continuous control in a planar motion subjected only to the Sun gravitational pulling can be formalized as
\begin{subequations}
	\begin{align}
		\Minimize_{\mathbf{u}(t)}\; J &:= -m(t_f)\\[0.6em]
		\shortintertext{subject to}
		\dot{\mathbf{x}}=\mathbf{f}\left(\mathbf{x},\mathbf{u}\right)=
		\begin{bmatrix}
			\dot{\mathbf{r}}\\
			\dot{\mathbf{v}}\\
			\dot{m}
		\end{bmatrix}&=
		\begin{bmatrix}
			\mathbf{v}\\
			-\mu\frac{\mathbf{r}}{\|\mathbf{r}\|^3}+\frac{T(t)}{m(t)}\boldsymbol{\gamma}\\
			-\frac{T(t)}{I_{\mathrm{sp}}g_0}
		\end{bmatrix} \label{eq:EoM1}\\
		\mathbf{x}(t_0)&=\mathbf{x}_0\\
		\mathbf{x}(t_f)&=\mathbf{x}_f \\
		0\le T(t)&\le T_{\max} \label{eq:u_cons}
	\end{align}
\end{subequations}

where $\mathbf{r}=\left[x, y\right]^T$, $\mathbf{v}=\left[u, v\right]^T$ denote the position and the velocity, respectively, $m$ the mass, $I_{\mathrm{sp}}$ represents the thruster specific impulse, $g_0$ gravitational acceleration at sea level, and the parameter $\mu$ is the Sun gravitational constant. The control $\mathbf{u}=\left[T, \alpha\right]^T$ is composed by the thrust magnitude $T$, having $T_{\max}$ as its attainable maximum value, and the angle $\alpha$, used to define the thrust direction $\boldsymbol{\gamma}=\left[\cos\alpha, \sin\alpha\right]^T$.\\
In order to facilitate the problem solution, the coupling between state and control in Eqs.~\eqref{eq:EoM1} is removed by defining a scaled control variable and a different metric for the mass \cite{wang2018minimum}

\begin{equation*}
	\Gamma(t)=\frac{T(t)}{m(t)}, \quad w(t)=\ln(m(t))
\end{equation*}

In this way, the ordinary differential equation for $m$ (i.e., the last term in Eq.~\eqref{eq:EoM1}) is substituted by
\begin{equation}
	\dot{w}=-\frac{\Gamma}{I_{\mathrm{sp}}g_0}
\end{equation}
while Eq.~\eqref{eq:u_cons} in the new variables is
\begin{equation}
	0\le\Gamma(t)\le T_{\max{}}e^{-w(t)}
	\label{eq:Gamma_nl}
\end{equation}

Following the SCP paradigm, the dynamics is linearized about a reference solution $(\bar{\mathbf{x}}, \bar{\mathbf{u}})$, and a trust-region constraint is added to keep the linearization valid \cite{malyuta2022convex}. Thus, the convexified optimization problem for a generic continuous-thrust problem reads
\begin{subequations}
	\begin{align}
		\Minimize_{\mathbf{u}(t)}\; J &:= -w(t_f)+\int_{t_0}^{t_f}\kappa\|\boldsymbol{\nu}(t)\|_1\,\mathrm{d}t+\int_{t_0}^{t_f}\theta \max\left(0, \eta(t)\right)\,\mathrm{d}t \label{eq:J_con}\\[0.6em]
		\shortintertext{subject to}
		\dot{\mathbf{x}}&=\mathbf{f}(\bar{\mathbf{x}}, \bar{\mathbf{u}}) +A(\mathbf{x}-\bar{\mathbf{x}})+B(\mathbf{u}-\bar{\mathbf{u}})+\boldsymbol{\nu}\label{eq:lin_dyn}\\
		\mathbf{x}(t_0)&=\mathbf{x}_0\\
		\mathbf{x}(t_f)&=\mathbf{x}_f \\
		0\le\Gamma&\le T_{\max{}}e^{-\bar{w}}\left(1-w+\bar{w}\right)+\eta
		\label{eq:Gamma_cons}\\
		\|\mathbf{x}-\bar{\mathbf{x}}\|_1+\tau\|\mathbf{u}-\bar{\mathbf{u}}\|_1&\le R
		\label{eq:linear_cons}
	\end{align}
\end{subequations}
where the dependence on time has been dropped for clarity's sake, and where
\begin{equation*}
	A=\frac{\partial \mathbf{f}}{\partial \mathbf{x}}\Bigg\rvert_{\bar{\mathbf{x}}}=
	\begin{bmatrix}
		\mathbf{0}_{2\times2} & I_2 & \mathbf{0}_{2\times1}\\
		-\mu\left(\frac{I_2}{\|\bar{\mathbf{r}}\|^3}-\frac{3\bar{\mathbf{r}}\bar{\mathbf{r}}^T}{\|\bar{\mathbf{r}}\|^5}\right) & \mathbf{0}_{2\times2} & \mathbf{0}_{2\times1}\\
		\mathbf{0}_{1\times2} & \mathbf{0}_{1\times2} & 0
	\end{bmatrix}
\end{equation*}
\begin{equation*}
	B=\frac{\partial \mathbf{f}}{\partial \mathbf{u}}\Bigg\rvert_{\bar{\mathbf{u}}}=
	\begin{bmatrix}
		\multicolumn{2}{c}{\textbf{0}_{2\times2}}\\
		\cos\bar\alpha & -\bar\Gamma\sin\bar\alpha\\
		\sin\bar\alpha & \bar\Gamma\cos\bar\alpha\\
		-1/{I_{\mathrm{sp}}g_0}&0
	\end{bmatrix}
\end{equation*}
with $I_n$ representing the $n$-dimensional identity matrix, and $\mathbf{0}_{n\times m}$ the $(n\times m)$ null array.\\
 Eq.~\eqref{eq:Gamma_cons} is the linearization about the reference solution of the non-convex constraints in Eq.~\eqref{eq:Gamma_nl}.\\
The variables $\boldsymbol{\nu}$ and $\eta$ are virtual controls introduced to avoid artificial infeasibility: $\boldsymbol{\nu}$ is added to the linearized dynamics constraints, while $\eta$ relaxes the linearized maximum-thrust constraint. Both terms are penalized in the objective function through the weights $\kappa$ and $\theta$, respectively, which are used to balance feasibility recovery against the original performance index. Artificial unboundeness is managed by imposing a constraint on the trust region radius $R$ (Eq.~\eqref{eq:linear_cons}).\\
The convex formulation of the optimal control problem must be translated in a programming problem to be solved. To this aim, the timespan $[t_0, t_f]$ is subdivided into $N$ segments, and both the state and the control are projected into this time grid to have finite-dimensional variables. A first-order hold (FOH) method is used to discretize the problem. Specifically, the control history is approximated as a piecewise affine function, that is in each segment \cite{malyuta2022convex}
\begin{equation}
	\mathbf{u}(t)=\frac{t_{k+1}-t}{h}\mathbf{u}_k+\frac{t-t_k}{h}\mathbf{u}_{k+1}=\lambda^-(t)\mathbf{u}_k+\lambda^+(t)\mathbf{u}_{k+1},\qquad \forall t\in\left[t_k,t_{k+1}\right]
\end{equation}
In this case, the linearized dynamics in Eq.~\eqref{eq:lin_dyn} becomes
\begin{equation}
	\dot{\mathbf{x}}=A\mathbf{x}+B\lambda^-\mathbf{u}_k+B\lambda^+\mathbf{u}_{k+1}+\left(\mathbf{f}\left(\bar{\mathbf{x}},\bar{\mathbf{u}}\right)-A\bar{\mathbf{x}}-B\bar{\mathbf{u}}\right)
	\label{eq:FOH_dyn}
\end{equation}
Hence, for each segment, Eq.~\eqref{eq:FOH_dyn} could be integrated by exploiting the linear systems theory. Thus, the dynamical constraint for each segment become \cite{malyuta2022convex}
\begin{equation}
	\mathbf{x}_{k+1}-\left(\hat{A}_k\mathbf{x}_k+\hat{B}_k^-\mathbf{u}_k+\hat{B}_k^+\mathbf{u}_{k+1}+\mathbf{q}_k+\boldsymbol{\nu}_k\right)=\mathbf{0}
	\label{eq:FOH_def}
\end{equation}
where
\begin{eqnarray}
	\hat{A}_k&=&\Phi(t_k, t_{k+1})\label{eq:linMatrices_0}\\
	\hat{B}_k^-&=&\hat{A}_k\int_{t_k}^{t_{k+1}}\Phi^{-1}(t_k, t)B\lambda^-\, \mathrm{d}t\\
	\hat{B}_k^+&=&\hat{A}_k\int_{t_k}^{t_{k+1}}\Phi^{-1}(t_k, t)B\lambda^+\, \mathrm{d}t\\
	\mathbf{q}_k&=&\hat{A}_k\int_{t_k}^{t_{k+1}}\Phi^{-1}(t_k, t)\left(\mathbf{f}\left(\bar{\mathbf{x}},\bar{\mathbf{u}}\right)-A\bar{\mathbf{x}}-B\bar{\mathbf{u}}\right)\, \mathrm{d}t
	\label{eq:linMatrices_end}
\end{eqnarray}
where $\Phi(t_k, t)$ is the state transition matrix between time $t_k$ and time $t$.\\
The optimal control problem can be translated so in a second-order cone programming problem (SOCP), i.e., a convex programming problem, that can be stated as:\\

\textbf{Problem 1 (SCP).\label{prob:1}} Find $\mathbf{y}=[\mathbf{x}_k, \mathbf{u}_k, \boldsymbol{\nu}_k, \eta_k], \ k=0,\dots,N$, such that
\begin{subequations}
	\begin{align}
		J :=& -w_N+\sum_{i=0}^N\kappa\|\boldsymbol{\nu}_k\|_1+\sum_{i=0}^N\theta\max\left(0, \eta_k\right)\label{eq:JJ}
		\shortintertext{is minimized, subjected to}
		\Delta&=\mathbf{x}_{k+1}-\left(\hat{A}_k\mathbf{x}_k+\hat{B}_k^-\mathbf{u}_k+\hat{B}_k^+\mathbf{u}_{k+1}+\mathbf{q}_k+\boldsymbol{\nu}_k\right)=\mathbf{0} &\forall& k\in1,\dots,N-1
		\label{eq:defects}\\
		0\le\Gamma_k&\le T_{\max{}}e^{-\bar{w}_k}\left(1-w_k+\bar{w}_k\right)+\eta_{k} &\forall& k\in1,\dots,N\label{eq:Tcon}\\
		&\|\mathbf{x}-\bar{\mathbf{x}}\|_1+\tau\|\mathbf{u}-\bar{\mathbf{u}}\|_1\le R\label{eq:RR}
	\end{align}
\end{subequations}
with $\mathbf{x}_0$ and $\mathbf{x}_N=\mathbf{x}_f$ prescribed.\\
The problem in Eqs.~\eqref{eq:JJ}--\eqref{eq:RR} must be solved iteratively. Iterations are stopped successfully if the maximum nonlinear constraint violation and the improvement in the cost function, i.e., $\|\mathbf{c}(\mathbf{x}, \mathbf{u})\|_\infty\le\varepsilon_c$ and $\bar{w}-w\le\varepsilon_w$. On the other hand, the procedure is terminated without success if the relative difference of the solution vector in two consecutive iterations is below a certain threshold (i.e., $\|\bar{x}-x\|_1\le\varepsilon_x$), meaning that the solution is not improving, or if the trust region radius is too small (i.e., $R\le\varepsilon_R$). As a matter of fact, following the SCP procedure, the trust radius must be changed according to an appropriate mechanism. At each iteration, the solution of the convex solver is accepted or rejected based on the following ratio \cite{malyuta2022convex}
\begin{equation}
	r=\frac{\bar\phi-\phi}{\bar\phi-\varphi}
	\label{eq:r}
\end{equation}
where $\phi$ is the actual cost, while $\varphi$ is the predicted cost reduction.
In this case, three parameters $(r_0, r_1, r_2)$, with $0<r_0<r_1<r_2<1$, are selected and used to manage the trust region. If $r\le r_0$, the  linearized problem is performing poorly, so the step is rejected and a smaller new region radius $R$ is selected. Otherwise, the solution is accepted and the trust region is updated as \cite{malyuta2022convex}
\begin{equation}
	R=\begin{cases}
		\begin{aligned}
			R/a \quad &\mbox{if} &r_0\le r<r_1\\
			R \quad  &\mbox{if} &r_1\le r<r_2\\
			b R \quad  &\mbox{if} &r\ge r_2
		\end{aligned}
	\end{cases}
	\label{eq:R}
\end{equation}
with $a$ and $b$ predetermined values grater than 1.

\subsection{Quantum Annealing Optimization}
Quantum annealing is a quantum computing technique developed to accelerate the solution of complex optimization problems.
QA frames a discrete optimization task as the search for the ground state of a quantum system whose energy landscape encodes the cost function to be minimized. Formally, an $n$-qubit register is prepared in the ground state of a driver Hamiltonian ${H}_{D}$ and evolves it slowly towards a problem Hamiltonian $H_{P}$ that embeds the objective function: 
\begin{equation}
	H(s)=\bigl(1-s\bigr)H_{D}+s\,H_{P},\qquad s=\tfrac{t}{\mathcal{T}}\in[0,1],
\end{equation}
with $\mathcal{T}$ the annealing time.  If the interpolation is sufficiently adiabatic, i.e., $\mathcal{T}\gg\mathcal{O}\left(\frac{1}{g_{\min}^{2}}\right)$ where $g_{\min}$ denotes the minimum spectral gap between the instantaneous ground and first-excited states, the system remains in its ground state, so that measuring the qubits at $s=1$ returns the optimal (or a near-optimal) solution with high probability \cite{farhi2000quantum}.\\
To interface with current hardware, the optimisation problem must be cast either as an Ising model, 
\begin{equation}
	J_{\text{Is}}(\mathbf{s})=\sum_{i}h_{i}s_{i}+\sum_{i<j}J_{ij}s_{i}s_{j},\qquad s_{i}\in\{-1,+1\}
\end{equation}
or, equivalently, as a quadratic unconstrained binary optimization (QUBO) problem, 
\begin{equation}
	J_{\text{QUBO}}(\mathbf{z})=\mathbf{z}^{\mathsf T}Q\,\mathbf{z},\qquad z_{i}\in\{0,1\}
\end{equation}

Both models are equivalent and it is possible to change from one model to the other via the simple linear transformation $z_{i}=(s_{i}+1)/2$. The resulting quadratic coefficients populate the entries of $H_{P}$ and therefore dictate the pairwise couplings implemented on the chip. Because physical qubits occupy a sparse connectivity graph, the logical QUBO must be embedded into the hardware topology using chains of ferromagnetically coupled qubits that behave as a single logical variable \cite{zbinden2020embedding}.  This embedding could inflate qubit counts and can introduce errors when chains break, so problem size and sparsity play a critical role in performance.\\
The possibility to encode several real-life optimization problems in QUBOs made the use of QA an appealing methodology to solve complex academic or industrial-scale instances \cite{yarkoni2022quantum}. Hence, in recent years, several advancements have been done to provide QA access to the community. To this aim, commercial annealers able to provide $\mathcal{O}(10^{3})$ qubits and millisecond-scale anneals, such as D-Wave's Advantage\footnote{\url{https://www.dwavequantum.com/solutions-and-products/systems/}. Last accessed on December 15, 2025.}, have been made available, enabling the solution of medium-scale optimization instances outright or, in hybrid modes, as sub-solvers inside larger classical algorithms. While gate-based quantum processors pursue universal speed-ups, QA already supplies a practical heuristic whose effectiveness hinges on judicious problem reformulation, careful scheduling, and clever  discretizations that keep QUBO density and embedding overhead within hardware limits.

\section{Low-Thrust Trajectory Quantum Optimization}\label{sec:QUBO}
To use a quantum annealer for low-thrust trajectory optimization, the optimal-control problem must be translated into a binary formulation that can be expressed as a QUBO. While the SCP paradigm is exploited to convert the continuous problem in a convex programming problem, as explained in Section \ref{subsec:SCP}, Problem 1 itself must be transformed in a QUBO form, trying to keep the number of variables contained in order to be compliant with the limitation in qubits given by the annealers.

\subsection{The Quadratic Unconstrained Binary Optimization Problem}
Problem 1 is a constrained problem, having both equality and inequality convex constraints. Thus, in order to get a QUBO, all the constrains must be encoded within the cost function, and at the same time being converted into quadratic forms. To have an unconstrained problem, a penalty method is exploited, meaning that the cost function $J$ is augmented with the constraints, i.e., $\hat{J} = J+\rho h$, where $\rho$ is a penalty coefficient and $h$ is a general equality constraint function, rearranged to be always positive and being null when the constraints are respected. As a matter of fact, when $\rho\rightarrow\infty$, the constraint $h$ is force to go to 0 and the original problem is correctly solved. Thus, all the constraints in Problem 1 must be modified to be compatible with this methodology.\\
In order to reduce the number of variables, the virtual buffer variables $\eta_k$ are dropped, and the artificial infeasibility is handled exploiting only the virtual control $\boldsymbol{\nu}$, as in \cite{kamath2023real}. Additionally, since the initial state and the final position and velocity are prescribed, the variables $\mathbf{x}_{1}$ and $\mathbf{x}_{N}$ are fixed to the desired values $\mathbf{x}_{0}$ and $\mathbf{x}_{f}$, respectively.\\
The equality constraint in Eq.~\eqref{eq:defects} could be easily transformed in a quadratic form. As a matter of fact, Eq.~\eqref{eq:defects} could be more conveniently rewritten as
\begin{equation}
	A_{eq}\mathbf{y}=\mathbf{b}_{eq}
\end{equation}
with
\begin{equation}
	A_{eq}=\left[\begin{array}{cccccc:cccccc:cccc}
		0 & -I & 0 & \dots & 0 & 0 & -\hat{B}_1^- & -\hat{B}_1^+ & 0 &\dots & 0 & 0 & -I & 0 & \dots & 0\\
		0 & \hat{A}_2 & -I & \dots & 0 & 0 & 0 & -\hat{B}_2^- & -\hat{B}_2^+ & \dots & 0 & 0 & 0 & -I & \dots & 0\\
		\vdots & \vdots & \vdots & \ddots & \vdots & \vdots & \vdots & \vdots & \vdots & \ddots & \vdots & \vdots & \vdots & \vdots & \ddots & \vdots\\
		0 & 0 & 0 & \dots & \hat{A}_{N-1} & -1 & 0 & 0 & 0 & \dots &  -\hat{B}_{N-1}^- & -\hat{B}_{N-1}^+ & 0 & 0 & \vdots & -I
	\end{array}\right]
	\label{eq:Aeq}
\end{equation}
and
\begin{equation}
	\mathbf{b}_{eq} = \left[\left(\mathbf{q}_1-\hat{A}_1\mathbf{x}_0\right)^T \quad \mathbf{q}_2^T \quad \dots \quad \left(\mathbf{q}_{N-1}+\mathbf{x}_f\right)^T\right]^T
	\label{eq:beq}
\end{equation}
Thus, the quadratic form to be considered in the augmented cost is simply
\begin{equation}
	J_{eq}=\left(A_{eq}\mathbf{y}-\mathbf{b}_{eq}\right)^T\left(A_{eq}\mathbf{y}-\mathbf{b}_{eq}\right)=\mathbf{y}^TA_{eq}^TA_{eq}\mathbf{y}-2\mathbf{b}_{eq}^TA_{eq}\mathbf{y}+\mathbf{b}_{eq}^T\mathbf{b}_{eq}
\end{equation}
Inequality constraints represented by Eq.~\eqref{eq:Tcon} must be first transformed into equality constraints, by means of slack variables, i.e.,
\begin{equation}
	G\mathbf{y}+\boldsymbol{\sigma}=\mathbf{h}
\end{equation}
where 
\begin{equation}
	G=\left[\begin{array}{ccccc:ccccc:cccc}
		0 & 0 & 0 &\dots & 0 & 1 & 0 & 0 & \dots & 0 & 0 &\dots & 0\\
		0 & \mathbf{g}_2 & 0 & \dots & 0 & 0 & 1 & 0 & \dots & 0 & 0 & \dots & 0\\
		0 & 0 & \mathbf{g}_3 & \dots & 0 & 0 & 0 & 1 & \dots & 0 & 0 & \dots & 0\\ 
		\vdots & \vdots & \vdots & \ddots & \vdots & \vdots & \vdots & \vdots & \ddots & \vdots & \vdots & \ddots & \vdots\\
		0 & 0 & 0 & \dots & T_{\max{}}e^{-\bar{z}_N} & 0 & 0 & 0 & \dots & 1 & 0 & \dots & 0
	\end{array}\right]
	\label{eq:G}
\end{equation}
and
\begin{equation}
	\mathbf{h}=-T_{\max{}}[\left(1+\bar{w}_1\right)e^{-\bar{w}_1}+T_{\max{}}e^{-\bar{w}_i}m_0 \quad \left(1+\bar{w}_2\right)e^{-\bar{w}_2} \quad \dots \quad \left(1+\bar{w}_N\right)e^{-\bar{w}_N}]^T
	\label{eq:h}
\end{equation}
with $\mathbf{g}_k = [0 \quad 0 \quad 0 \quad 0 \quad T_{\max{}}e^{-\bar{w}_k}]$ and $\boldsymbol{\sigma}$ the slack variables which must be positive to keep the new equality constraint equivalent to the original inequality one. Thus, the quadratic form to be considered in the augmented cost is simply
\begin{equation}
	J_{ineq}=\left(G\mathbf{y}+\boldsymbol{\sigma}-\mathbf{h}\right)^T\left(G\mathbf{y}+\boldsymbol{\sigma}-\mathbf{h}\right)=\mathbf{y}^TG^TG\mathbf{y}+\boldsymbol{\sigma}^T\boldsymbol{\sigma}+2\boldsymbol{\sigma}^TG\mathbf{y}-2\mathbf{h}^TG\mathbf{y}-2\mathbf{h}^T\boldsymbol{\sigma}+\mathbf{h}^T\mathbf{h}
\end{equation}
Lastly, the trust region constraint in Eq.~\eqref{eq:RR} must be adapted to the augmented cost. While the use of the $\ell^1$-norm, as standard in the SCP for low-thrust transfer, is still possible to this aim, it would require to define additional slack variables, that would inflate the number of qubits. For this reason, a different approach is preferred, that is to define the trust region using the $\ell^2$-norm, since it can be directly translated in a quadratic form. Thus, Eq.~\eqref{eq:RR} is first rewritten as
\begin{equation}
	\|\mathbf{x}-\bar{\mathbf{x}}\|_2+\tau\|\mathbf{u}-\bar{\mathbf{u}}\|_2\le \frac{R}{\sqrt{n_x+n_u}}
\end{equation}
where $n_x$ and $n_u$ are the length of the vector $\mathbf{x}$ and $\mathbf{u}$, respectively. The radius $R$ is scaled to make the comparison between the standard SCP and the QUBO easier. Then, it is converted into the quadratic form
\begin{equation}
	J_R = \mathbf{y}^TG_R\mathbf{y}+\mathbf{g}_R^T\mathbf{y}+h_R+\sigma_R
\end{equation}
where 
\begin{equation}
	\renewcommand{\arraystretch}{0.8}
	G_R=\begin{bmatrix}
		I_{n_x} & &\\
		& \tau^2I_{n_u} &\\
		& & 0
	\end{bmatrix}
	\label{eq:GR}
\end{equation}
\begin{equation}
	\mathbf{g}_R=\left[-2\bar{\mathbf{x}}^T \quad -2\tau^2\bar{\mathbf{u}}^T \quad \mathbf{0}_{n_\nu}^T\right]^T
	\label{eq:gR}
\end{equation}
and $h_R=\bar{\mathbf{x}}^T\bar{\mathbf{x}}+\tau^2\bar{\mathbf{u}}^T\bar{\mathbf{u}}-R^2/(n_x+n_u)$, with $\sigma_R$ a slack variable, constrained to be always positive.\\
Finally, the cost in Problem 1 (i.e., Eq.~\eqref{eq:JJ}) must be rearranged. As for the trust region constraint, while the use of $\ell^1$-norm and max functions could be still possible, it will require additional variables. For this reason, they are both changed into the squared $\ell^2$-norm, that is
\begin{equation}
	J=-z_N+\kappa^2\boldsymbol{\nu}^T\boldsymbol{\nu}=\mathbf{y}^TQ\mathbf{y}+\mathbf{q}^T\mathbf{y}
\end{equation}
with
\begin{equation*}
	\renewcommand{\arraystretch}{0.7}
	Q=\begin{bmatrix}
		0 & &\\
		& 0 &\\
		& & \kappa^2I_{n_\nu}
	\end{bmatrix}
\end{equation*}
and $\mathbf{q}=\left[\mathbf{0}_{n_x-1}^T \quad -1 \quad \mathbf{0}_{n_u}^T \quad \mathbf{0}_{n_\nu}^T\right]^T$.\\
In this way, the convex programming problem is translated in a quadratic unconstrained problem (QUP), that can be stated as:\\
\textbf{Problem 2 (QUP).} Find $\mathbf{y}=[\mathbf{x}, \mathbf{u}, \boldsymbol{\nu}]$, such that
\begin{equation}
	\begin{aligned}
		J:=\mathbf{y}^TQ\mathbf{y}+\mathbf{q}^T\mathbf{y}\ +\ &\rho_1\left(\mathbf{y}^TA_{eq}^TA_{eq}\mathbf{y}-2\mathbf{b}_{eq}^TA_{eq}\mathbf{y}+\mathbf{b}_{eq}^T\mathbf{b}_{eq}\right)+\\
		&\rho_2\left(\mathbf{y}^TG^TG\mathbf{y}+\boldsymbol{\sigma}^T\boldsymbol{\sigma}+2\boldsymbol{\sigma}^TG\mathbf{y}-2\mathbf{h}^TG\mathbf{y}-2\mathbf{h}^T\boldsymbol{\sigma}+\mathbf{h}^T\mathbf{h}\right)+\\
		&\rho_3\left(\mathbf{y}^TG_R\mathbf{y}+\mathbf{g}_R^T\mathbf{y}+\sigma_R+h_R\right)
	\end{aligned}
	\label{eq:J2}
\end{equation} 
is minimized, with $\Gamma_k\ge0$, $\boldsymbol{\sigma}\ge0$, and $\sigma_R\ge0$.\\
The weighting parameters $\rho_1$, $\rho_2$, and $\rho_3$ are selected such that there is a good balance between the cost functions and the constraints and among the different constraints.\\
To derive the QUBO representation, each continuous decision variable must be encoded into $n_b$ binary variables using a fixed-point scheme that maps the real-valued domain onto the allowable binary range. Specifically, the minimum and maximum representable values of the chosen binary encoding are selected to be aligned with the lower and upper bounds of each variable. These bounds come either from the optimization problem (e.g., for $\Gamma$) or from physical considerations. If $y_{i, \min}$ and $y_{i,\max}$ are the real limits, the fixed-point mapping ensures that the binary vector $\mathbf{z}_i\in\{0,1\}^{n_b}$
corresponds to
\begin{equation}
	y_i = y_{i, \min} +\frac{y_{i,\max}-y_{i, \min}}{2^{n_b}-1}\,\sum_{k=0}^{n_b-1}2^k z_{i,k}
	\label{eq:binenc}
\end{equation}
such that $\mathbf{z}_i=\mathbf{0}$ yields $y_{i, \min}$ and $\mathbf{z}_i=\mathbf{1}$ yields $y_{i, \max}$. This procedure converts each real variable into a fixed-length binary vector whose span exactly covers the feasible interval, enabling substitution into the objective and constraint expressions. In a more compact way, starting from this definition, it is possible to define a linear vector map between the continuous real variable $\mathbf{y}$ and the binary variable vector $\mathbf{z}$
\begin{equation}
	\mathbf{y} = \mathbf{y}_{\min}+D\mathbf{z}
	\label{eq:ztoy}
\end{equation}
where
\begin{equation}
	D = \frac{\mathbf{y}_{\max}-\mathbf{y}_{\min}}{2^{n_b}-1}
	\begin{bmatrix}
		\mathbf{d}^T & 0 & \dots & 0\\
		0 & \mathbf{d}^T & \dots & 0\\
		\vdots & \vdots & \ddots & \vdots\\
		0 & 0 & \dots & \mathbf{d}^T
	\end{bmatrix}
	\label{eq:D}
\end{equation}
with $\mathbf{d}=\left[2^0 \quad 2^1 \quad \dots \quad 2^{n_b-1}\right]$. The same transformation can be applied also to the slack variables to get
\begin{equation}
	\boldsymbol{\sigma} = \boldsymbol{\sigma}_{\min}+D_\sigma\mathbf{z}_\sigma, \qquad \sigma_R=\sigma_{R,\min}+D_{\sigma_R}\mathbf{z}_{\sigma_R}
	\label{eq:ztos}
\end{equation}
By inserting Eqs.~\eqref{eq:ztoy} and \eqref{eq:ztos} in Eq.~\eqref{eq:J2}, it is possible to define an expression for the cost function in terms of the binary decision vector $\hat{\mathbf{z}}=\left[\mathbf{z}^T \quad \mathbf{z}_\sigma^T \quad \mathbf{z}_{\sigma_R}^T\right]^T\in\{0, 1\}$, that is
\begin{equation}
	J:=\hat{\mathbf{z}}^T\hat{Q}\hat{\mathbf{z}}+\hat{\mathbf{q}}^T\hat{\mathbf{z}}+\hat{k}
\end{equation}
where $\hat{Q}$ is the matrix, composition of the quadratic terms, $\hat{\mathbf{q}}$ is the sum of the linear terms, and $k$ is the known coefficient. Since for binary variables $z_i^2=z_i$, it is possible to incorporate the linear terms within the quadratic matrix, setting $\tilde{Q}=\hat{Q}+\diag\left(\hat{q}\right)$. The constant coefficient $k$, on the other hand, can be ignored in the optimization setting.\\
Hence, it is possible to render the quadratic unconstrained problem as a quadratic unconstrained binary problem, that can be finally stated as:\\
\textbf{Problem 3 (QUBO).} Find $\hat{\mathbf{z}}\in\{0, 1\}$, such that
\begin{equation}
	J:=\hat{\mathbf{z}}^T\tilde{Q}\hat{\mathbf{z}}
\end{equation}
is minimized.

\subsection{The Computational Algorithm}
As per Problem 1, Problem 3 must be solved iteratively to obtain the optimal solution. However, some care must be placed to have the quantum-based sequential convex programming (quSCP) to converge to the non-linear optimal solution.\\
In classical constrained optimization, quadratic penalties enforce feasibility by solving a sequence of unconstrained problems. Convergence to the constrained optimum relies on letting the penalty coefficient starting with small values and progressively increase them, enforcing stricter constraints as the algorithm progresses. In this continuous setting, violations can be made arbitrarily small and the optimal solution is retrieved when $\rho\to\infty$. In the case of quantum optimization, the decision variables are discretized via a binary expansion (as per Eq.~\eqref{eq:binenc}). The optimization is then performed over a finite set of binary configurations, and the goal is not to approximate the constrained optimum by increasing the penalty coefficient, but to construct a single unconstrained binary energy whose ground state corresponds to the optimal feasible discretized solution. Thus, the role of the penalty coefficients is fundamentally different, since it is a separation parameter used to ensure that all infeasible binary states have higher energy than any feasible state, while preserving objective discrimination among feasible states.\\
Discretization changes constraint enforcement because constraint residuals become quantized. Over the discrete set, each residual component $\left(A_{eq}\mathbf{y}-\mathbf{b}_{eq}\right)_j$ (and likewise $\left(G\mathbf{y}+\boldsymbol{\sigma}-\mathbf{h}\right)_j$) can only take values on a finite grid determined by the minimum jump for each variable $\Delta_i=\frac{y_{i,\max}-y_{i, \min}}{2^{n_b}-1}$ and the constraint coefficients. As a result, there typically exists a strictly positive \emph{minimum nonzero violation} $\delta$ provided the discretization and scaling. This gap is the key property that enables an exact penalty separation in the discrete setting. Unlike continuous optimization, infeasibility cannot approach zero continuously, but it must jump from zero to at least $\delta$.\\
A sufficient strategy to choose the initial $\rho$ values is to select its order of magnitude so that the smallest possible penalty incurred by any infeasible configuration dominates the largest possible improvement in the original objective across the discretized domain. In this case, any
\begin{equation*}
	\rho>\frac{\Delta w_N}{\delta^2}
\end{equation*} 
with $\Delta w_N$ the maximum possible variation in the objective function, guarantees that every global minimizer of $J_{\text{QUBO}}$ 
is feasible, and among feasible points minimizing $J_{\text{QUBO}}$ is equivalent to minimizing the original cost function, i.e., the control effort. In practice, $\Delta w_N$ is selected to be conservative, by computing the difference between the initial mass and the final mass with the thrust always on. Additionally, a safety margin $\zeta>1$ is introduced, so that
\begin{equation}
	\rho = \zeta\frac{\Delta w_N}{\delta^2}
	\label{eq:rho}
\end{equation}
to account both for mathematical issues and for hardware characteristics, e.g., analog noise, coefficient quantization, and embedding-related distortion.\\
It is important to remark that Eq.~\eqref{eq:rho} provides a sufficient separation condition for the ideal discrete QUBO, but it should not be interpreted as a prescription to select arbitrarily large penalty values on quantum hardware. Current annealers operate within a finite programmable range for local biases and couplers; therefore, the QUBO matrix must be globally rescaled before submission to the solver. If the penalty coefficients are chosen excessively large, this rescaling compresses the contribution of the original objective function, reducing the effective energy separation among feasible states. As a consequence, although infeasible configurations are energetically discouraged, the propellant-minimization term may become comparable to, or smaller than, the effective resolution imposed by analog noise, coefficient quantization, and embedding-related distortions. This creates a trade-off between feasibility enforcement and objective resolution. In this work, the penalty factors are therefore selected to be sufficiently large to separate feasible and infeasible configurations in the discretized model, but not increased beyond what is required for reliable feasibility, so as to preserve the relative contribution of the original cost after hardware rescaling.\\
While the fact that QUBO formulation requires the discretization of the decision variables via an unsigned binary expansion can be helpful to a certain extent in the selection of proper penalty factors, it generates issues in the representation of the optimal solution. As a matter of fact, even with an appropriate selection of the penalty factors, the solution of the QUBO in Problem 3 will never match the solution of Problem 2, and, consequently, the one of the original Problem 1. For a fixed binary depth $n_b$, the discretization step $\Delta_i$ imposes an intrinsic quantization error that limits the accuracy with which the continuous optimum can be represented. In practice, conservative initial bounds $\left[\mathbf{y}_\text{min}, \mathbf{y}_\text{max}\right]$ are often required to ensure that the feasible region is not inadvertently excluded. However, wide bounds increase $\Delta_i$ and therefore degrade solution fidelity. This trade-off is particularly detrimental when the continuous optimum lies in a small subset of the initial hyper-rectangle, since a large portion of the binary code space is spent representing values that are never selected by optimal solutions. Iterative refinement addresses this limitation by repeatedly re-centering and shrinking the encoding interval around a high-quality QUBO solution, thereby reducing $\Delta_i$ without increasing the representing bits. Once the QUBO has been solved with the initial encoding, a candidate continuous vector $\mathbf{y}^\star$
is obtained by decoding the best (i.e., lowest-energy feasible) binary sample. The encoding interval is then shrunk by a factor $\beta\in\left(0, 1\right)$. Defining the current width of each variable as $W_i = y_{i,\max}-y_{i, \min}$ the refined width is
\begin{equation}
	W_i=\max\left(\beta\, W_i, W_{i, \min}\right)
	\label{eq:Wi}
\end{equation}
with $W_{i, \min}$ a prescribed value to prevent over-refinement and the new discretization step
\begin{equation}
	\Delta_i=\frac{W_i}{2^{n_b}-1}
\end{equation}
To ensure that the optimal value is exactly representable by the unsigned encoding at the next iteration, a correcting index $\gamma=\lfloor\frac{2^{n_b}-1}{2}\rfloor$ is introduced and the boundaries are computed as
\begin{equation}
	y_{i, \min}=y_i^\star-\gamma\Delta_i, \qquad y_{i,\max}=y_{i, \min}+\left(2^{n_b}-1\right)\Delta_i
	\label{eq:bnds}
\end{equation}
This construction places the optimal values at the center of the grid, while decreasing the discretization step by the factor $\beta$. 
If global physical bounds $\left(L_i, U_i\right)$ are available, the interval is clipped to satisfy $y_{i, \min}\ge L_i$ and $y_{i,\max}\le U_i$, adjusting $\gamma$ if necessary to maintain feasibility of the interval. Because interval refinement changes both the attainable objective range over the discrete domain and the quantization of constraint residuals, the penalty weight used to enforce feasibility in the QUBO should be recomputed at each iteration to preserve feasibility--optimality separation. This iterative refinement procedure yields a sequence of progressively finer QUBO encodings with fixed bit-depth, enabling increasingly accurate approximations of the continuous constrained optimum.\\
It should be noted that this refinement strategy is local by construction. Since the encoding interval is progressively re-centered and shrunk around the best decoded QUBO solution, an incorrect early sample may bias the subsequent search and cause the refinement to converge toward a suboptimal discrete basin. This risk is particularly relevant when the QUBO solver returns approximate samples, as is generally the case for quantum annealing and hybrid heuristics. To mitigate this effect, the refinement is not performed around an arbitrary sample, but around the lowest-energy feasible solution returned by a batch of samples. Moreover, the shrink factor $\beta$ is selected conservatively, so that the refined interval still preserves a neighborhood around the current solution rather than collapsing immediately to a narrow region. If the decoded objective does not improve, if the best solution repeatedly lies close to the boundary of the current interval, or if no feasible sample is recovered after refinement, the shrink step is rejected and the interval is either restored to the previous bounds or enlarged before repeating the QUBO solve. In this sense, the refinement should be interpreted as a trust-region-like local improvement mechanism rather than as a guarantee of global optimality.\\
The iterative refinement could be performed, however, only if the solution of the QUBO could be accepted, that is when the trust region constraint, needed to guarantee the validity of the dynamics linearization, is respected. Since the trust-region constraint is embedded within the cost function, an approach similar to the soft trust region strategy used in SCP is exploited. Specifically, the violation with respect to the trust-region constraint is evaluated, i.e., $\delta_R=\left(\mathbf{y}_{\min}+D\mathbf{z}\right)^TG_R\left(\mathbf{y}_{\min}+D\mathbf{z}\right)+\mathbf{g}_R^T\left(\mathbf{y}_{\min}+D\mathbf{z}\right)-h_R$. If the constraint is violated (i.e., $\delta_R>0$, thus the step has been outside the trust region), the step is rejected and the penalty coefficient is increased, setting $\rho_3=\zeta_r\rho_3$. Otherwise, the step is accepted and the a new QUBO, with updated bounds, as per Eqs.~\eqref{eq:bnds}, is solved.\\
The process is repeated until a stopping criterion is met, i.e., $\max{\Delta_i}$ falls below a prescribed tolerance, the decoded objective ceases to improve, or the step between two consecutive optimal value is too small. Once the solution of the single linearized problem represented by Problem 3 has been found, the validity of the solution for the single SCP iteration is evaluated and, in case, the new reference solution is updated and the update rule for the radius in Eq.~\eqref{eq:R} is followed. This procedure is, again, repeated iteratively, until the stopping criteria, which are the same of Problem 1 in Section \ref{subsec:SCP}, are met.\\
In conclusion, quSCP is made by two loops:
\begin{itemize}
	\item An \textbf{inner loop} solving the QUBO in Problem 3 several times, exploiting the iterative refinement to obtain a discrete solution close to the continuous one. This loop is equivalent to the iterative cycles performed by any convex solver;
	\item An \textbf{outer loop} solving the SCP, by updating the reference solution, adjusting the trust region, and checking the converge conditions. This loop is standard to any SCP methodology. 
\end{itemize}
A sketch summarizing how quSCP works is reported in Algorithm \ref{alg:quSCP}.

\begin{algorithm}[htp!]
	
	\caption{Quantum-based Sequential Convex Programming (quSCP)}
	\label{alg:quSCP}
	
	\setstretch{1}
	\DontPrintSemicolon
	\SetAlgoLined
	\LinesNumbered
	\SetArgSty{textnormal}
	
	\SetKwFor{For}{for}{}{end}
	\SetKwFor{While}{while}{}{end}
	\SetKwIF{If}{ElseIf}{Else}{if}{}{else if}{else}{end}
	\SetKwComment{Comment}{$\triangleright$}{}
	\SetKw{Break}{break}
	
	\KwIn{Spacecraft data, $\mathbf{x}_0$, $\mathbf{x}_f$, $t_0$, $t_f$, SCP parameters, $N$, $n_b$}
	Compute an initial guess for $\bar{\mathbf{x}}$ and $\bar{\mathbf{u}}$\;
	Set the initial value for the key penalties $\rho$\Comment*{See Eq.~\eqref{eq:rho}}
	Set the initial values for the boundaries $y_{i, \min}$ and $y_{i, \max}$\;
	\While{\text{\texttt{SCP\_exitconditions} is not \textbf{true}}}{
		Compute the linearized matrices for the dynamics $\hat{A}_k$, $\hat{B}_k^+$, $\hat{B}_k^-$, and $\mathbf{q}_k$\Comment*{See Eqs.~\eqref{eq:linMatrices_0}--\eqref{eq:linMatrices_end}}
		Build $A_{eq}$ and $\mathbf{b}_{eq}$\Comment*{See Eqs.~\eqref{eq:Aeq}--\eqref{eq:beq}}
		Build $G$ and $\mathbf{h}$\Comment*{See Eqs.~\eqref{eq:G}--\eqref{eq:h}}
		Build $G_R$, $\mathbf{g}_R$, and $h_R$\Comment*{See Eqs.~\eqref{eq:GR}--\eqref{eq:gR}}
		\While{\text{\texttt{QUBO\_exitconditions} is not \textbf{true}}
		}{
			Compute the enconding matrices $D$, $D_\sigma$, and $D_{\sigma_R}$\Comment*{See Eq.~\eqref{eq:D}}
			Define $\tilde{Q}$\;
			Solve the QUBO in Problem 3\;
			Retrieve $\mathbf{y}^\star$, $\boldsymbol{\sigma}^\star$, and $\boldsymbol{\sigma}_R^\star$\Comment*{See Eqs.~\eqref{eq:ztoy}--\eqref{eq:ztos}}
			\eIf(\Comment*[f]{Trust region respected}){$\delta_R\le0$}{
				Recompute the variable widths $W_i$ and new boundaries $y_{i, \min}$ and $y_{i, \max}$\Comment*{See Eqs.~\eqref{eq:Wi}--\eqref{eq:bnds}}
			}(\Comment*[f]{Trust region not respected}){
				Update the thrust region coefficient $\rho_3=\zeta_r\rho_3$\;
			}
			Compute the \texttt{QUBO\_exitconditions}\;
		}
		Compute the ratio $r$\Comment*{See Eq.~\eqref{eq:r}}
		\eIf(\Comment*[f]{SCP solution is rejected}){$r<r_0$}{
			Reduce trust region radius with $R=R/a$\;
		}(\Comment*[f]{SCP solution is accepted}){
			Update reference values $\bar{\mathbf{x}}$ and $\bar{\mathbf{u}}$ based on $\mathbf{y}^\star$\;
			Update trust region radius based on Eq.~\eqref{eq:R}\;
			Update the new boundaries $y_{i, \min}$ and $y_{i, \max}$\;
		}
		Compute the \texttt{SCP\_exitconditions}\;
	}
	\KwOut{Optimal trajectory $\mathbf{x}^\star$ and optimal control ${\mathbf{u}}^\star$}
\end{algorithm}

\section{Results}\label{sec:Results}
\begin{figure}[!htpb]
	\centering
	
	\subfigure[Optimal trajectory.\label{subfig:SCP100_traj}]{
		\includegraphics[width=0.5\linewidth]{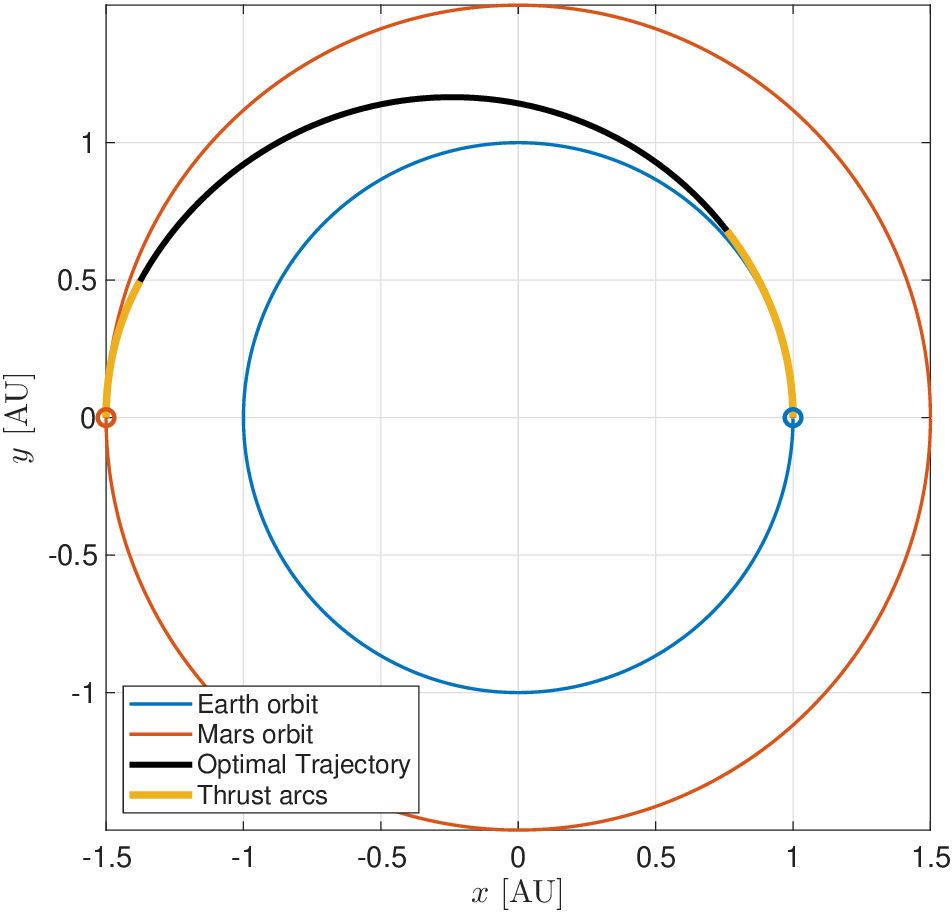}
	}
	\hfill
	\subfigure[Optimal control profile.\label{subfig:SCP100_ctrl}]{
		\includegraphics[width=0.5\linewidth]{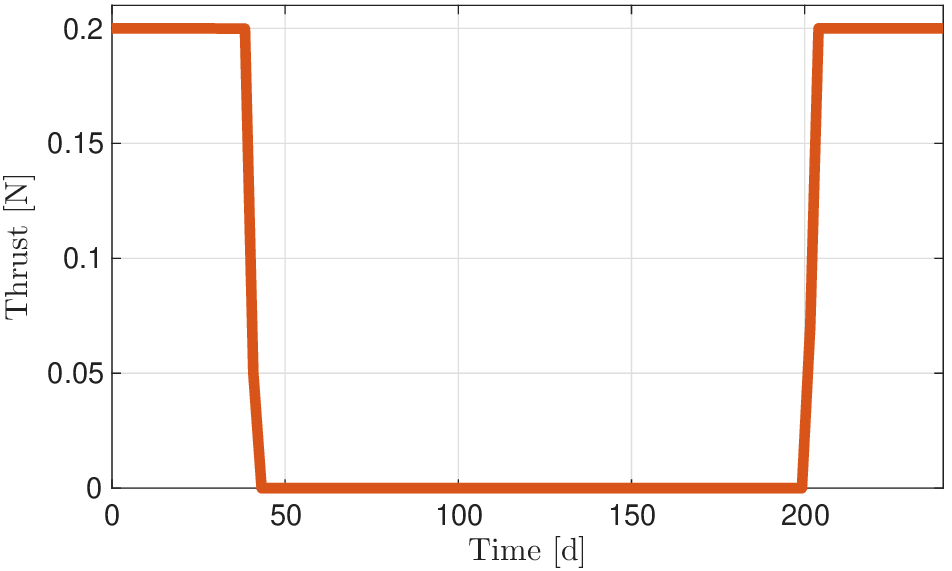}
	}
	\caption{SCP optimal solution with $N=100$.}
	\label{fig:SCP100}
\end{figure}

The methodology was validated on an optimal Earth--Mars transfer problem. Characteristics of the initial and final state and of the spacecraft are listed in Table \ref{tab:scenariodata}. The resulting solutions are compared against a reference trajectory computed with a standard SCP solver, whose parameters are listed in Table \ref{tab:SCPpar}.\\
To characterize the capabilities and performance of quSCP, its key internal parameters, namely the discretization size $N$ and the bit-width $n_b$ are varied. The QUBO subproblems generated within quSCP were solved using D-Wave's cloud-accessible optimization backends, leveraging both hardware quantum annealing solvers (QPU-based) and hybrid solvers. In particular, QPU-based solutions, which require a minor-embedding of the logical QUBO onto the physical connectivity graph, are best suited to relatively small instances, while hybrid solvers, which combine classical metaheuristics with targeted QPU calls are exploited to handle larger problem sizes. The distinction is practically relevant because, as the number of logical variables or the QUBO density increases, the embedding typically induces longer chains and may ultimately become infeasible, whereas hybrid approaches mitigate these limitations and retain robust performance on higher-dimensional instances. For QPU runs, logical problems were embedded onto the hardware graph using Ocean's standard minor-embedding workflow; chain strength and unembedding strategy were selected consistently across experiments, and embedding quality was monitored via chain-length and chain-break statistics. All QUBOs were globally rescaled to comply with solver coefficient ranges while preserving the relative weighting between the original objective and penalty terms. Solutions were selected as the lowest-energy feasible samples (or, when multiple feasible samples were returned, the one minimizing the original objective), using a fixed sampling budget (i.e., number of reads) for fair comparison across configurations. All solver calls, both for QPU and hybrid ones, were orchestrated through the D-Wave Ocean SDK, which provides a unified interface for model construction, embedding, sampling, and post-processing.\\
The numerical results are not intended to establish a runtime advantage over classical SCP or QUP solvers. Instead, they are used to assess whether the proposed QUBO transcription can produce feasible and physically consistent low-thrust trajectories, to quantify the degradation introduced by binary discretization and current quantum hardware, and to identify the regimes in which QPU-based and hybrid solvers can be practically employed.\\
The benchmark is provided by an SCP solution with $N=100$, shown in Fig.~\ref{fig:SCP100}.\\

\begin{table}[ht!]
	\caption{Data for the test scenario.}
	\label{tab:scenariodata}
	\centering
	\begin{tabular}{l c c}
		\hline
		\textbf{Parameter} & \textbf{Symbol} & \textbf{Value} \\
		\hline
		Spacecraft mass & $m$ & \qty{250}{\kilogram}\\
		Maximum thrust & $T_{\max}$ & \qty{0.2}{\newton}\\
		Specific impulse & $I_\textbf{sp}$ & \qty{1500}{\second}\\
		Initial position & $\mathbf{r}_0$ & \qtylist{1.496e8; 0.0}{\kilo\meter}\\
		Initial velocity & $\mathbf{v}_0$ & \qtylist{0.0; 29.78}{\kilo\meter\per\second}\\
		Final position & $\mathbf{r}_f$ & \qtylist{-2.244e8; 0.0}{\kilo\meter}\\
		Final velocity & $\mathbf{v}_f$ & \qtylist{0.0; -24.32}{\kilo\meter\per\second}\\
		Time of flight & $t_f-t_0$& \qty{240}{\day}\\
		\hline
	\end{tabular}
\end{table}

\begin{table}[ht!]
	\caption{SCP parameters.}
	\label{tab:SCPpar}
	\centering
	\begin{tabular}{c r}
		\hline
		\textbf{Parameter} & \textbf{Value} \\
		\hline
		$\kappa$ & 1000\\
		$\tau$ & 100\\
		$R_0$ & 2500\\
		$\left(a, b\right)$ & $\left(1.5, 1.5\right)$\\
		$\left(r_0, r_1, r_2\right)$ & $\left(0.01, 0.5, 0.95\right)$\\
		$\left(\varepsilon_c, \varepsilon_w, \varepsilon_x, \varepsilon_R\right)$ & $\left(10^{-10}, 10^{-4}, 10^{-8}, 10^{-4}\right)$\\
		\hline
	\end{tabular}
\end{table}

\begin{figure}[!htpb]
	\centering
	\includegraphics[width=0.5\linewidth]{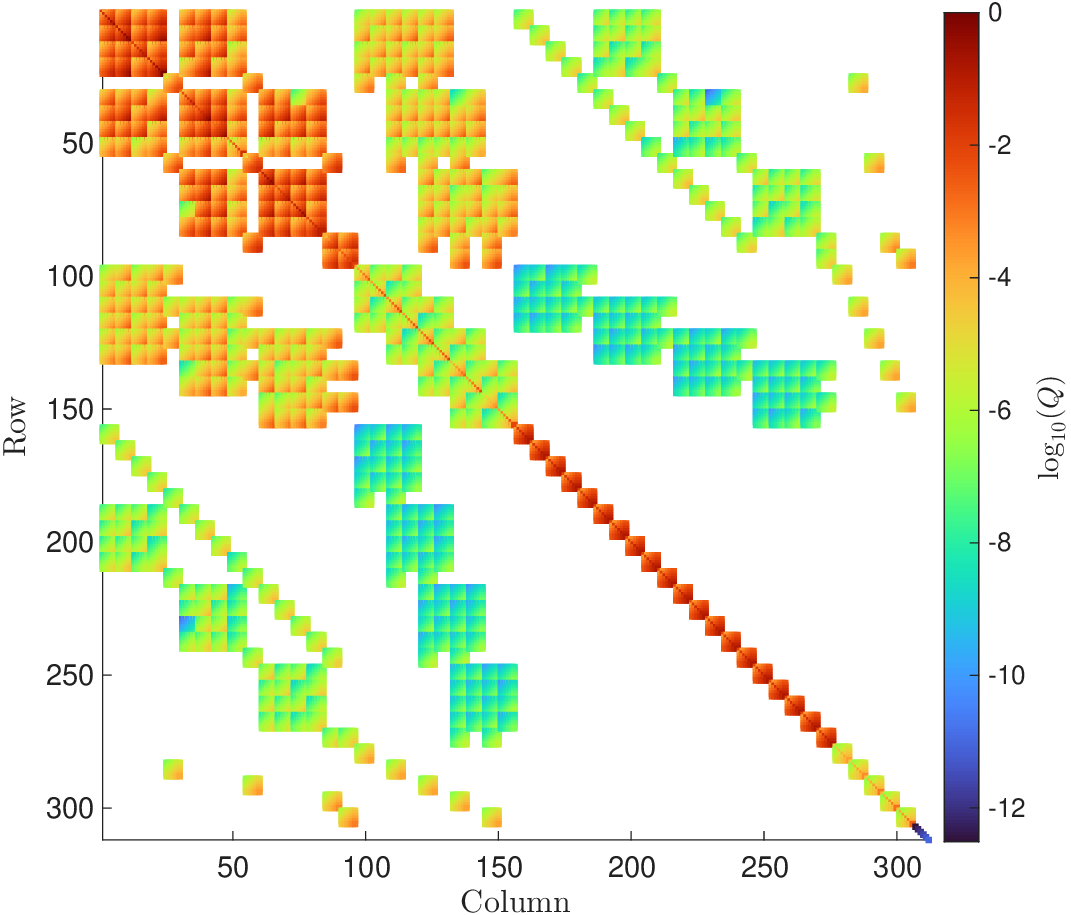}
	\caption{Structure for the normalized matrix $\tilde{Q}$ of the QUBO problem with $N=4$ and $n_b=6$.}
	\label{fig:QN4np6}
\end{figure}

\begin{figure}[!htpb]
	\centering
	\includegraphics[width=0.75\linewidth]{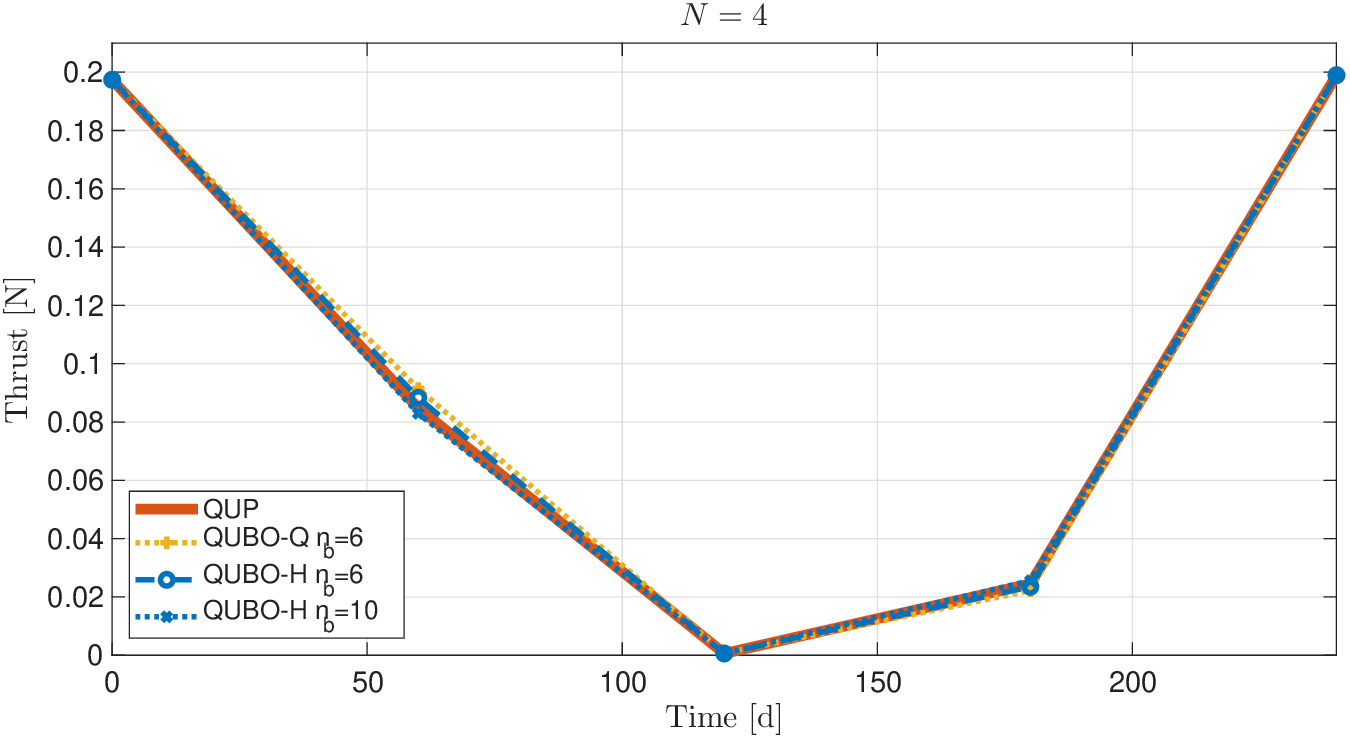}
	\caption{Optimal control profile for different solvers with $N=4$.}
	\label{fig:Control4}
\end{figure}

As anticipated, QPU-based sampling imposes strict limitations on problem size: increasing dimensionality leads to longer minor-embedding chains, which amplifies the impact of analog noise and chain breaks and can in turn degrade solution quality. On the other hand, a low bit-width $n_b$ results in a coarse discretization of the decision variables, which fundamentally limits the achievable accuracy of the solution by introducing significant quantization error, requiring a longer and more complex iterative refinement, while using a low number of segments $N$ could bring to poor control and dynamics representation.\\
Indeed, for $N=4$ and $n_b=6$, the resulting optimization problem comprises 312 binary variables. Figure \ref{fig:QN4np6} illustrates an example of the corresponding QUBO matrix $\tilde{Q}$, with nonzero entries highlighted using a color scale proportional to their magnitude. As shown, the quadratic formulation is relatively sparse, exhibiting an overall sparsity of approximately 0.24, and can be efficiently embedded and solved using a QPU-based sampler. Figure \ref{fig:Control4} depicts the optimal control profile obtained by quSCP with the QUBO formulation on the QPU (quSCP-Q), alongside the reference solution from the original SCP (with SOCP), the solution of the SCP with continuous QUP, solved with Clarabel, and the hybrid sampler result (quSCP-H). It can be seen that the history in time of the thrust magnitude does not change significantly for the different solvers. However, due to the low number of segments considered, it is quite different from the benchmark profile in Fig.~\ref{subfig:SCP100_ctrl}.\\
Qualitative considerations are reflected in the quantitative metrics reported in Table~\ref{tab:results4}. For $N=4$, all approaches yield comparable values of the maximum constraint violation $\lVert\mathbf{c}\rVert_\infty$, indicating that the nonlinear constraints are satisfied to a similar accuracy across solvers. As expected, the SCP benchmark attains the lowest propellant consumption, representing the continuous reference solution. The QUP formulation, solved with a classical convex optimizer, exhibits a modest increase in propellant mass, which can be attributed to the unconstrained formulation. The QUBO-based formulations incur a further, yet limited, degradation: for $n_b=6$, the QPU-based solution differs from the SCP benchmark by approximately 5.2\% in propellant mass, while maintaining a comparable level of constraint satisfaction. The hybrid solver achieves slightly improved propellant values with respect to the pure QPU approach and benefits from increased bit-width, as evidenced by the reduction in $m_p$ when moving from $n_b=6$ to $n_b=10$. This trend confirms the impact of discretization resolution on solution quality. From a computational standpoint, the classical SCP and QUP solvers complete in fractions of a second, whereas quSCP requires a substantially larger wall time, dominated by the repeated solution of QUBO subproblems on remote quantum hardware. Notably, the QPU-based configuration exhibits a significantly lower wall time per iteration, owing to the intrinsic speed of quantum annealing. However, it also requires a larger number of iterations to converge, reflecting the reduced quality of each QUBO solve, which is affected by long embedding chains, analog noise, and chain-break phenomena. In contrast, the hybrid approach trades longer per-iteration times for improved robustness and solution quality, resulting in fewer iterations overall. These results highlight the complementary nature of QPU and hybrid solvers within quSCP, i.e., while QPU sampling enables rapid exploration at small scales, hybrid strategies provide enhanced reliability and accuracy as problem complexity increases.

\begin{table}[ht!]
	\caption{Summary of results.}
	\label{tab:results4}
	\centering
	\begin{tabular}{l c c c c c c}
		\hline\hline
		$N$ & \textbf{Solver} & $n_b$ & $m_p$ [kg] & $\left\|\mathbf{c}\right\|_\infty$ [-]
		& \textbf{\# Iterations} & \textbf{WT} [s]\\
		\hline
		\multirow{5}{*}{4} 
		& SCP (SOCP)  & - & 106.23 & 0.0217 & -   & 0.125\\
		& SCP (QUP)  & - & 111.07 & 0.0218 & 48  & 0.259\\
		& quSCP-Q  & 6 & 111.75 & 0.0225 & 201 & 64.22\\
		& \multirow{2}{*}{quSCP-H} 
		& 6 & 111.57 & 0.0221 & 105 & 109.2\\
		& & 10 &111.15 & 0.0220 & 102 & 144.7\\
		\hline\hline
	\end{tabular}
\end{table}

\begin{figure}[!htpb]
	\centering

	\subfigure[$N=6$\label{subfig:Control6}]{
		\includegraphics[width=0.48\linewidth]{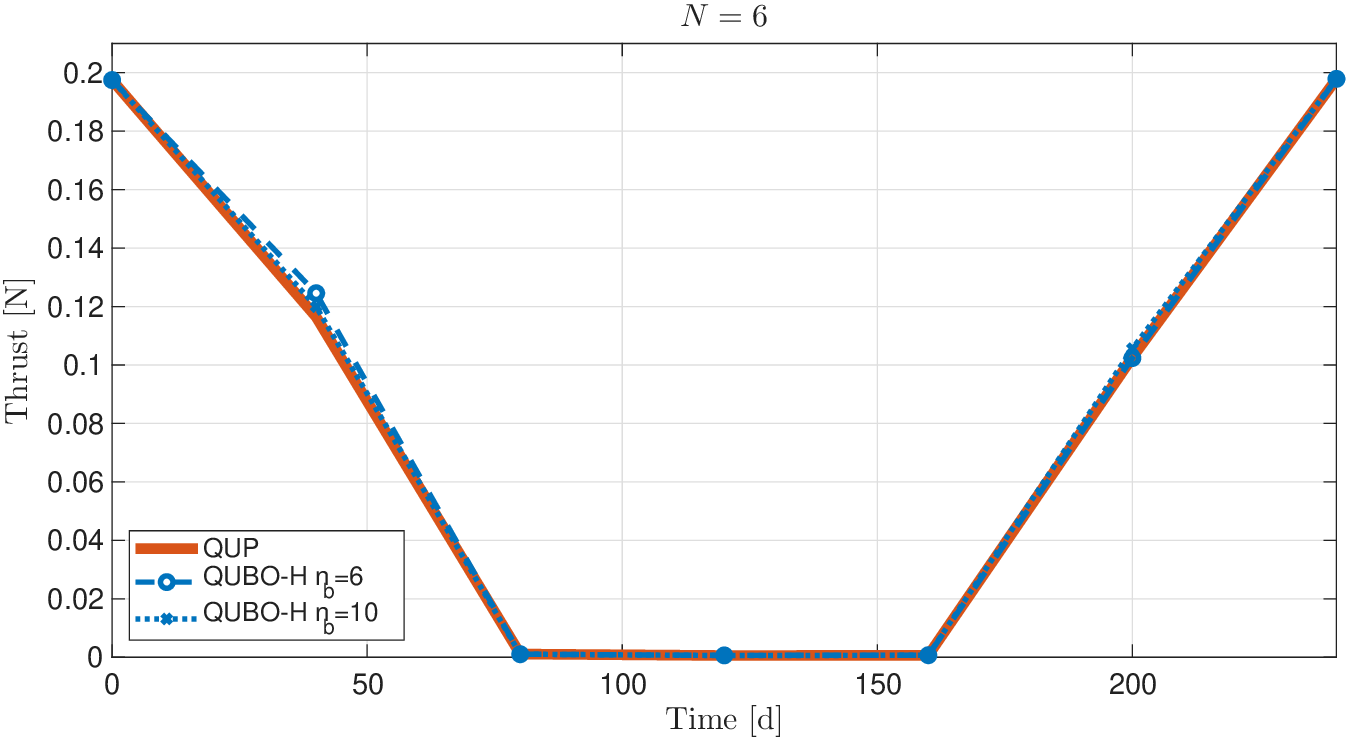}
	}
	
	\subfigure[$N=8$\label{subfig:Control8}]{
		\includegraphics[width=0.48\linewidth]{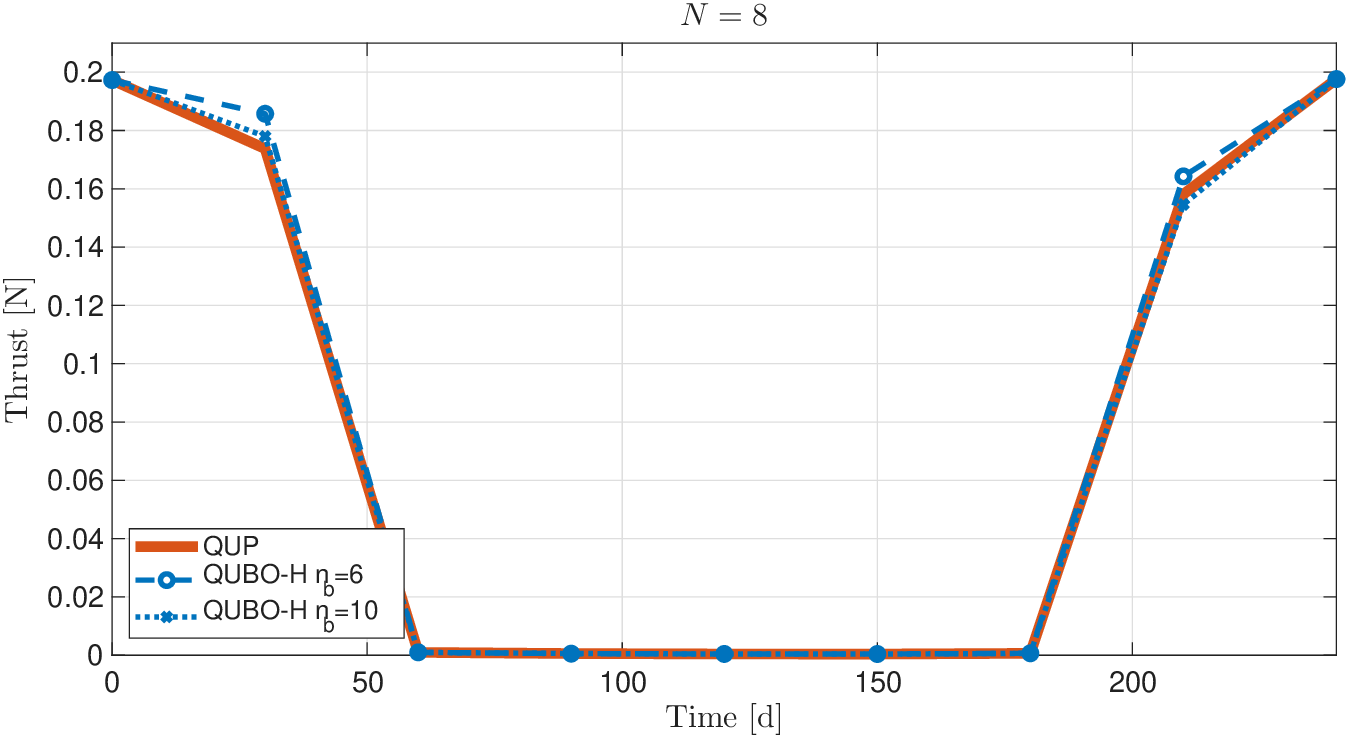}
	}
	\hfill
	\subfigure[$N=10$\label{subfig:Control10}]{
		\includegraphics[width=0.48\linewidth]{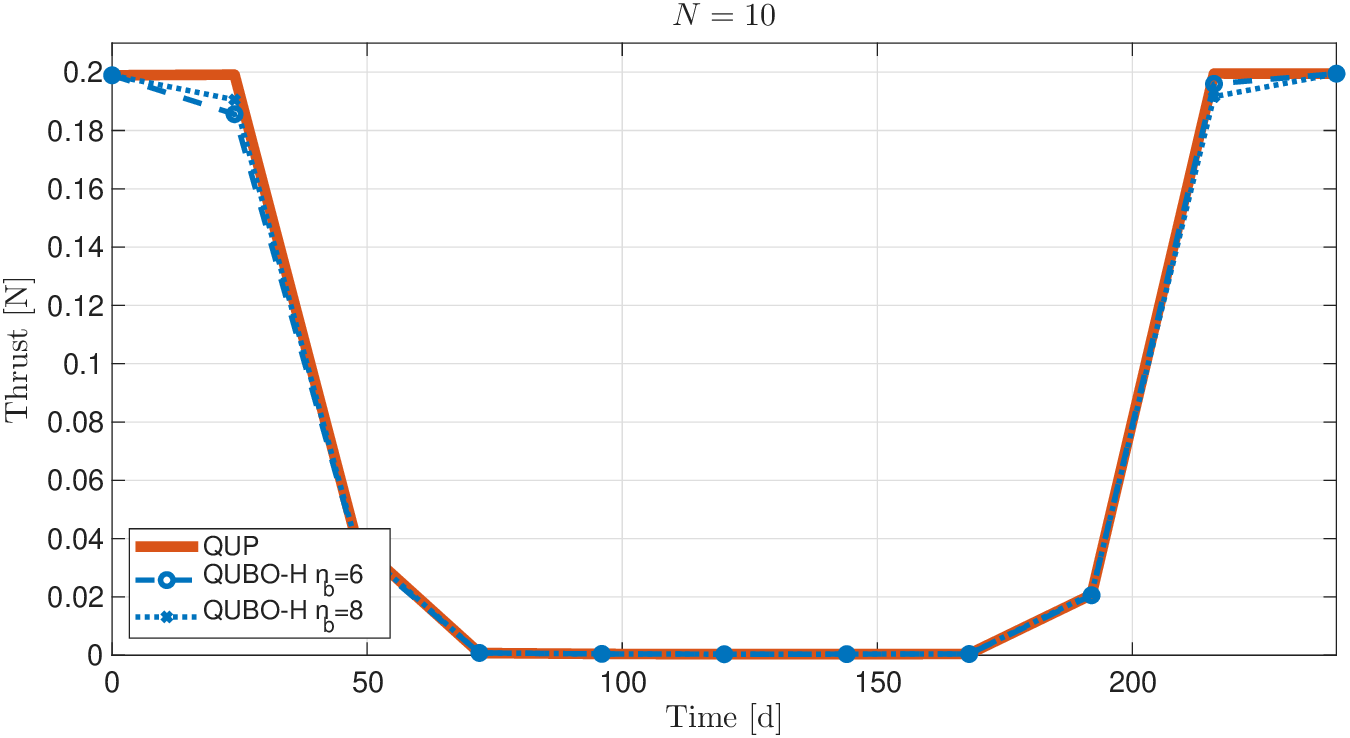}
	}
	
	\caption{Optimal control profile for different number of segments.}
	\label{fig:ControlAll}
\end{figure}

The same Earth--Mars transfer has been tested by increasing the number of segments. In these cases, the dimensionality of the problem grows accordingly and, as a consequence, it becomes no longer possible to embed the resulting QUBO on the QPU, thereby precluding the use of a QPU-based sampler. On the other hand, it is still possible to exploit hybrid solvers to prove that they are able to tackle this case study problem and offer solutions similar to SCP. The optimal control history for the hybrid-based samplers, together with SCP with SOCP and QUP formulation, for a number of segments from 6 to 10 are reported in Fig.~\ref{fig:ControlAll}. Quantitative results are summarized in Table~\ref{tab:results} for $N=6,8,10$. As the number of segments increases, the SCP benchmark exhibits the expected monotonic reduction in propellant mass and constraint violation, approaching the finely discretized reference at $N=100$. The QUP formulation closely tracks this behavior, with marginal differences in both $m_p$ and $\left\|\mathbf{c}\right\|_\infty$, confirming that the convex relaxation preserves high fidelity with respect to the SCP solution. The quSCP solved with hybrid sampler remains consistently aligned with the classical solvers across all tested values of $N$. For each discretization level, the hybrid solutions exhibit only a modest degradation in propellant mass with respect to SCP, on the order of 1\%-2\% for $N=6$ and below 1\% for $N=8$ and $N=10$, while maintaining comparable levels of constraint satisfaction. Increasing the bit-width systematically improves solution quality. For instance, at $N=6$, moving from $n_b=6$ to  $n_b=10$ reduces the propellant gap with SCP by approximately \qty{0.3}{\kilogram}, a trend that is consistently observed at larger $N$. From a computational standpoint, the classical SCP and QUP solvers complete in sub-second times, whereas quSCP requires substantially longer wall times, dominated by the iterative solution of QUBO subproblems. The hybrid solver converges in roughly 100 to 140 iterations, with wall times ranging from approximately \qtyrange{140}{310}{\second}, as $N$ increases. Despite this overhead, the hybrid framework demonstrates robust scalability and delivers solutions that closely approximate the classical optimum, even in regimes where direct QPU-based sampling is no longer feasible. These results confirm that hybrid quantum-classical solvers provide a viable path for extending quSCP to higher-dimensional trajectory optimization problems while retaining solution quality comparable to state-of-the-art classical methods.

\begin{table}[ht!]
	\caption{Summary of results.}
	\label{tab:results}
	\centering
	\begin{tabular}{l c c c c c c}
		\hline
		$N$ & \textbf{Solver} & $n_b$ & $m_p$ [kg] & $\left\|\mathbf{c}\right\|_\infty$ [-]
		& \textbf{\# Iterations} & \textbf{WT} [s]\\
		\hline
		\multirow{4}{*}{6} 
		& SCP (SOCP) & - & 96.77 & $1.706\times10^{-4}$ & -  & 0.098\\
		& SCP (QUP)  & - & 99.30 & $1.774\times10^{-4}$ & 57 & 0.315\\
		& \multirow{2}{*}{quSCP-H} 
		& 6 & 99.98 & $1.854\times10^{-4}$ & 135 & 240.2\\
		& & 10 & 99.67 & $1.801\times10^{-4}$ & 118 & 140.8\\
		\hline
		\multirow{4}{*}{8} 
		& SCP (SOCP) & - & 94.13 & $3.075\times10^{-6}$ & -  & 0.117\\
		& SCP (QUP)  & - & 94.25 & $4.138\times10^{-6}$ & 32 & 0.427\\
		& \multirow{2}{*}{quSCP-H} 
		& 6 & 94.86 & $5.052\times10^{-6}$ & 122 & 256.2\\
		& & 10 & 94.36 & $4.301\times10^{-6}$ & 97  & 175.4\\
		\hline
		\multirow{4}{*}{10} 
		& SCP (SOCP) & - & 91.87 & $6.165\times10^{-7}$ & -  & 0.128 \\
		& SCP (QUP) & - & 91.97 & $7.364\times10^{-7}$ & 52 & 0.578\\
		& \multirow{2}{*}{quSCP-H} 
		& 6 & 92.54 & $7.995\times10^{-7}$ & 136 & 313.6\\
		&     & 8 & 92.11 & $7.784\times10^{-7}$ & 128 & 208.5\\
		\hline
		100 & SCP (SOCP)  & - & 91.84 & $7.370\times10^{-9}$ & - &2.582\\
		\hline
	\end{tabular}
\end{table}

\section{Conclusions}\label{sec:Conclusions}
This work has introduced quSCP, a novel framework that integrates sequential convex programming with QUBO-based optimization to enable the solution of low-thrust trajectory design problems using quantum and hybrid quantum-classical solvers. Each convex subproblem generated within SCP is reformulated as a QUBO and solved through quantum-oriented backends, while the trust-region logic and convergence properties of SCP are preserved. The methodology has been validated on an Earth--Mars low-thrust transfer, demonstrating that quSCP produces physically consistent trajectories and closely tracks the solutions obtained with classical SCP and QUP formulations. Both QPU-based and hybrid solvers can be effectively embedded within the iterative pipeline, with hybrid strategies proving robust for higher-dimensional instances where direct QPU sampling is no longer feasible.\\
The study also exposes the current technological limitations of quantum hardware. The size of QUBO instances that can be directly embedded on present-day QPUs remains severely constrained, and solution quality is affected by embedding overhead, long chains, and analog noise. Consequently, the present implementation does not provide a computational advantage over state-of-the-art classical solvers, either in terms of accuracy or wall-clock time. This limitation is expected under current hardware and software conditions: direct QPU sampling is restricted to small QUBO instances because of embedding constraints, while hybrid solvers improve scalability at the cost of substantially larger runtimes. The contribution of this work is therefore methodological rather than computational-performance driven. As a matter of fact, it demonstrates that the SCP subproblems arising in low-thrust trajectory optimization can be consistently embedded into a QUBO framework and solved by quantum-oriented backends, producing solutions that remain close to classical references. These results indicate that, although current hardware does not yet enable a practical quantum advantage, quSCP provides a future-proof architecture for trajectory optimization. It establishes a concrete pathway for seamlessly integrating more capable quantum devices as they become available, positioning quantum computing as a viable component of next-generation mission design workflows.

\appendix
\subsection*{Acknowledgements}
\label{acknowledgements}
This work was supported by CASTOR, a project that has received funding from the European Union's Horizon Europe research and innovation programme under the Marie Skłodowska-Curie grant agreement no.~101103826.\\
The author acknowledges also the CINECA award under the ISCRA initiative, for the availability of quantum computing resources and support.

\subsection*{Declaration of competing interest}

The author has no competing interests to declare that are relevant to the content of this article.

\section*{References}

\bibliographystyle{astrobib}
\bibliography{bibliography}

\subsection*{Author biography}
\begin{biography}[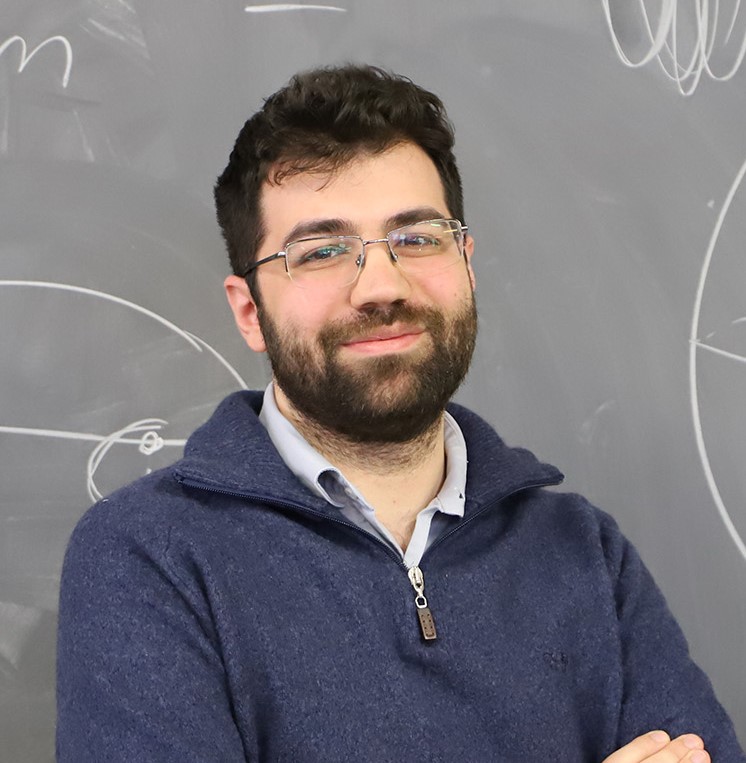]{Carmine Giordano} Carmine Giordano is an Assistant Professor at Politecnico di Milano. His main fields of interest are nonlinear astrodynamics and trajectory optimization. He took part in different space mission studies funded by ESA (LUMIO, Hera’s Milani, and RAMSES' Farinella) as mission analyst.
\end{biography}

\vspace*{2.6em}
\subsection*{Graphical table of contents}

\begin{figure*}[h!]
	\centering
	
	\subfigure[$N=6$]{
		\includegraphics[width=0.48\linewidth]{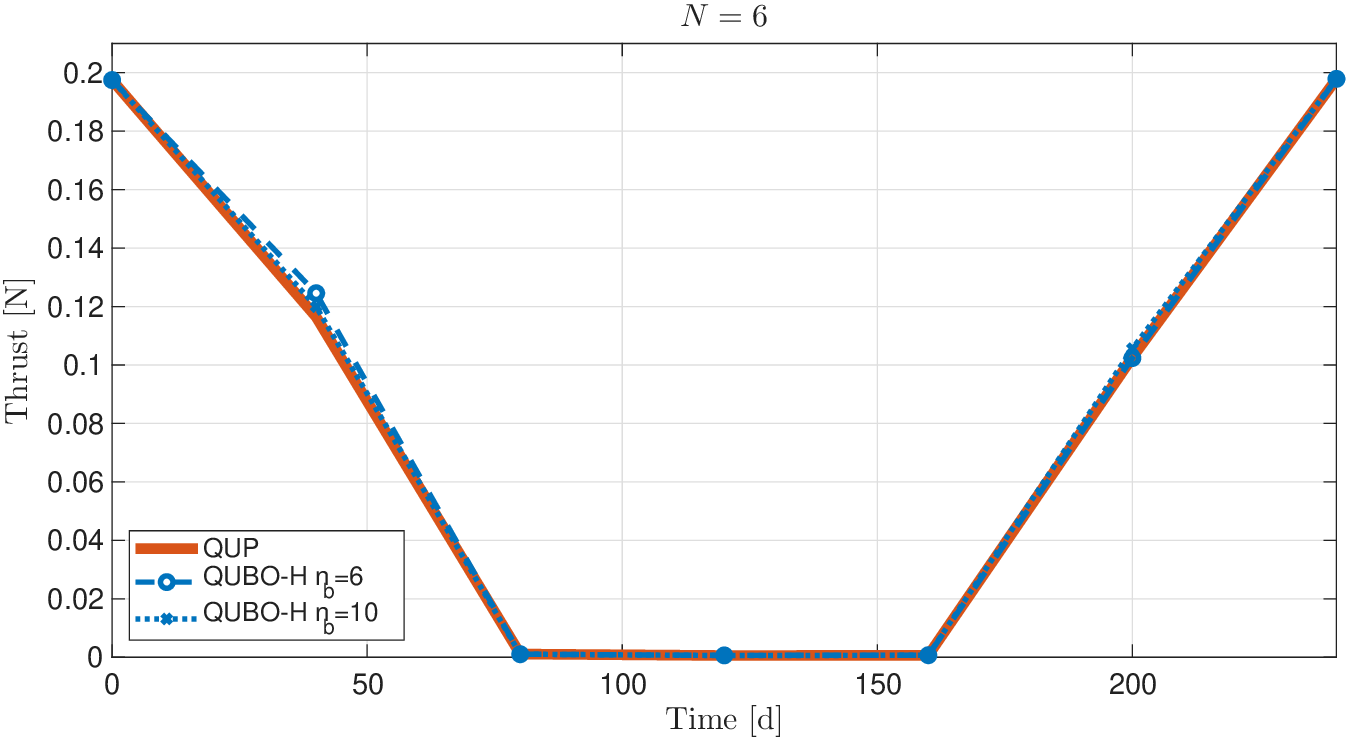}
	}
	
	\subfigure[$N=8$]{
		\includegraphics[width=0.48\linewidth]{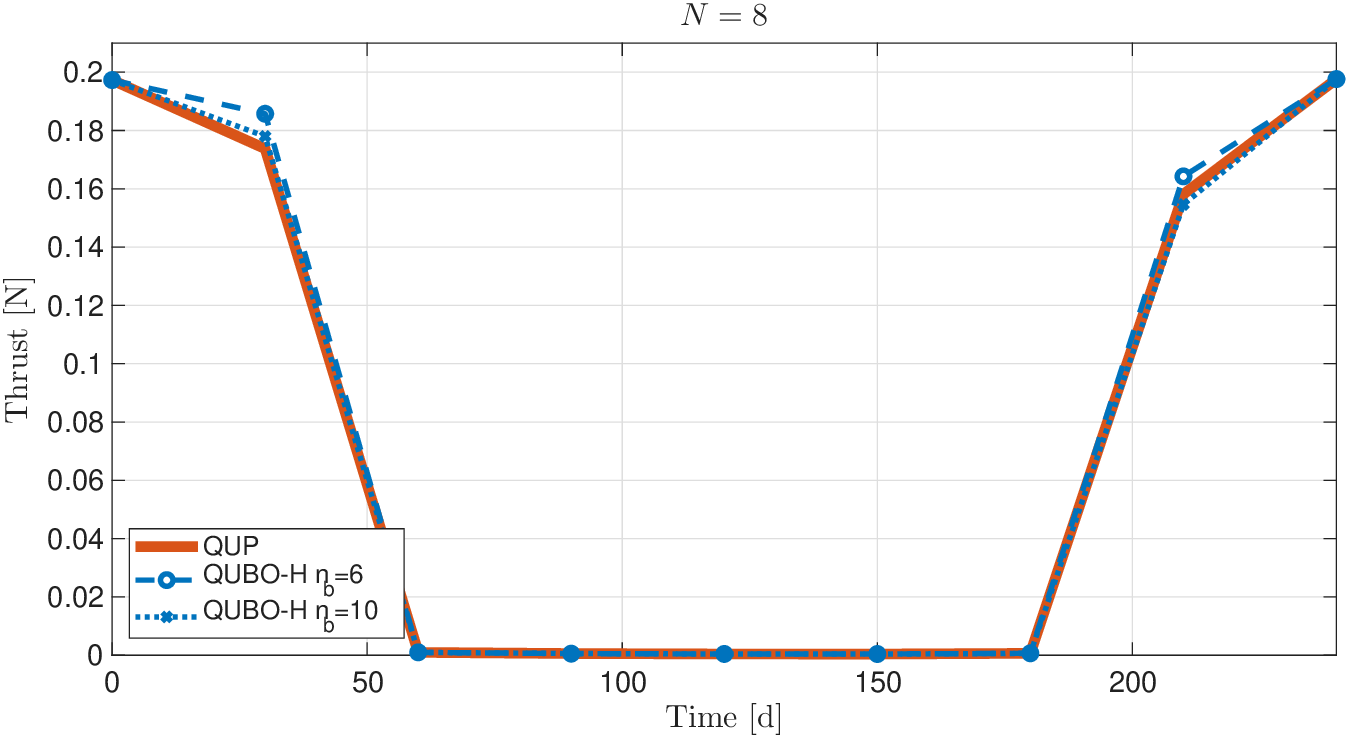}
	}
	\hfill
	\subfigure[$N=10$]{
		\includegraphics[width=0.48\linewidth]{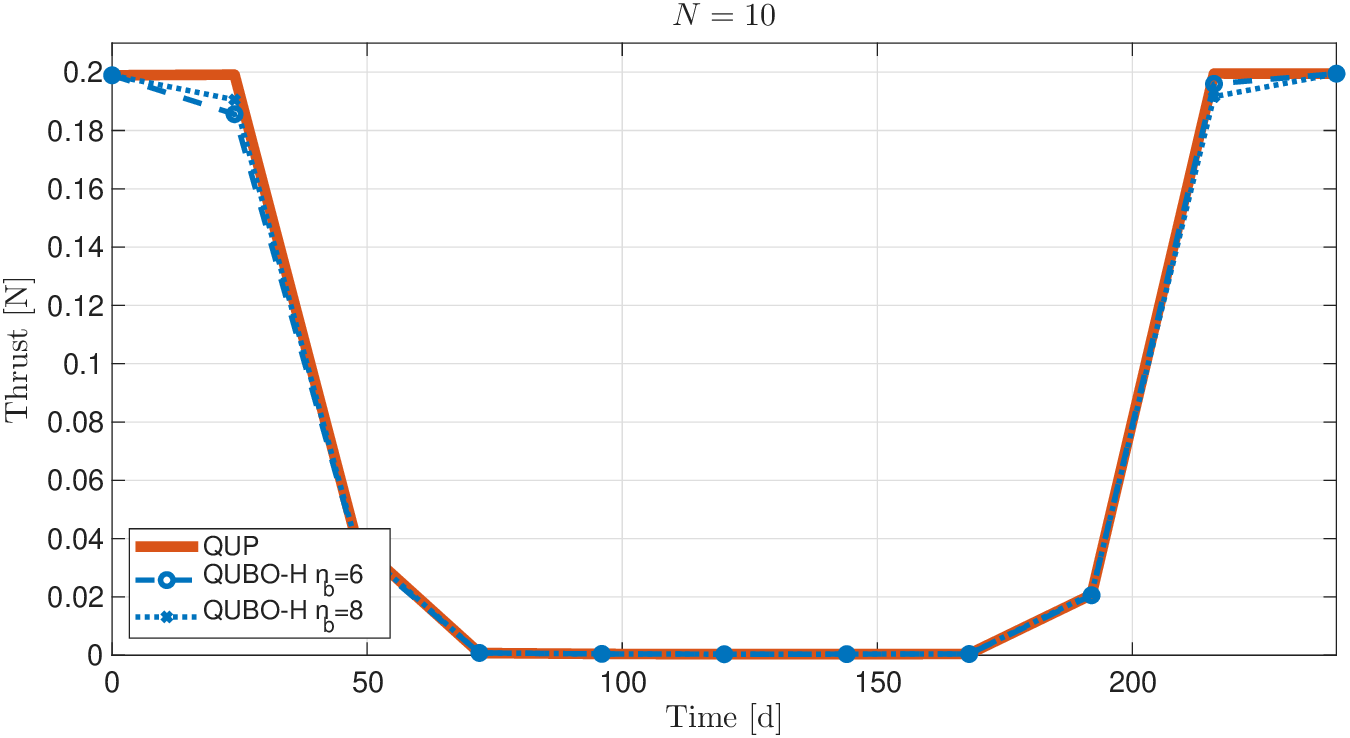}
	}
	
	\caption*{The figure illustrates the optimal low-thrust control profiles obtained with the proposed quSCP framework for increasing trajectory discretization levels. The close agreement with the classical QUP reference highlights the capability of hybrid quantum--classical solvers to recover physically consistent solutions for interplanetary trajectories.}
\end{figure*}

\end{document}